\documentclass[12pt]{article}

\usepackage[T1]{fontenc}
\usepackage[dvips]{graphicx}
\usepackage[arrow, matrix, curve]{xy}
\usepackage{amssymb}
\newtheorem{lemma}{Lemma}[section]
\newtheorem{theorem}[lemma]{Theorem}
\newtheorem{corollary}[lemma]{Corollary}
\newtheorem{prop}[lemma]{Proposition}

\newcommand{\ol}{\overline}

\def\Ind#1#2{#1\setbox0=\hbox{$#1x$}\kern\wd0\hbox to 0pt{\hss$#1\mid$\hss}
\lower.9\ht0\hbox to 0pt{\hss$#1\smile$\hss}\kern\wd0}
\def\Notind#1#2{#1\setbox0=\hbox{$#1x$}\kern\wd0\hbox to 0pt{\mathchardef
\nn="3236\hss$#1\nn$\kern1.4\wd0\hss}\hbox to 0pt{\hss$#1\mid$\hss}\lower.9\ht0
\hbox to 0pt{\hss$#1\smile$\hss}\kern\wd0}
\def\ind{\mathop{\mathpalette\Ind{}}}

\begin{document}
\title{Neostability-properties of Fra\"{\i}ss\'e limits of 2-nilpotent groups of exponent $p> 2$}
\author{Andreas Baudisch}
\date{\today}
\maketitle

\begin{abstract}
\noindent Let $L(n)$ be the language of group theory with $n$ additional new constant symbols $c_1,\ldots,c_n$. In $L(n)$ we consider the class ${\mathbb K}(n)$ of all finite groups $G$ of exponent $p > 2$, where  $G'\subseteq\langle c_1^G,\ldots,c_n^G\rangle \subseteq Z(G)$ and $c_1^G,\ldots,c_n^G$ are linearly independent. Using amalgamation we show the existence of   Fra\"{\i}ss\'e limits $D(n)$ of  ${\mathbb K}(n)$. $D(1)$ is Felgner's extra special $p$-group. The elementary theories of the $D(n)$ are superstable of SU-rank 1. They have the independence property.
\end{abstract}

\section{Introduction}

We consider the variety ${\mathbb G}_{2,p}$ of nilpotent groups of class 2 of exponent $p>2$ in the language $L$ of group theory. To get the Amalgamation Property (AP) in \cite{Bau} an additional predicate $P(G)$ for $G\in{\mathbb G}_{2,p}$ with $G'\subseteq P(G)\subseteq Z(G)$ is  introduced. Let ${\mathbb G}_{2,p}^P$ be the category of this groups in the extended language $L_P$ where the morphisms are embeddings. Using the class ${\mathbb K}_{2,p}^P$ of finite structures in ${\mathbb G}_{2,p}^P$ we get a Fra\"{\i}ss\'e limit $D$. If we build $D$ by amalgamation then $P(a)$ says that $a$ will become an element of the commutator subgroup $D'$ of $D$ in that process. In \cite{Bau} it is shown that ${\rm Th}(D)$ is not simple. Here we point out that $D$ has the tree property of the second kind (TP$_2$). This is easily seen.
\medskip

Let $L(n)$ be the language of group theory with $n$ additional new constant symbols $c_1,\ldots,c_n$. In $L(n)$ we consider the class ${\mathbb G}(n)$ of all groups $G\in {\mathbb G}_{2,p}$, where  $G'\subseteq\langle c_1^G,\ldots,c_n^G\rangle \subseteq Z(G)$ and $c_1^G,\ldots,c_n^G$ are linearly independent. We use linear independence, since we can consider an abelian group of expoent p as a vector space over ${\mathbb F}_p$.  $\langle X \rangle$ denotes the substructure generated by $X$. Hence   $\langle c_1^G,\ldots,c_n^G\rangle = \langle \emptyset \rangle$. ${\mathbb G}(n)$ is uniformly locally finite. Let ${\mathbb K}(n)$ be the class of finite structures in ${\mathbb G}(n)$. ${\mathbb K}(n)$ has the Hereditary Property (HP), the Joint Embedding Property (JEP) and the Amalgamation Property (AP). Hence the Fra\"{\i}ss\'e limit $D(n)$ of the class ${\mathbb K}(n)$ exists. Note that $D(1)$ is the extra special $p$-group considered by U.~Felgner in \cite{Fe}. In \cite{MacSt} the corresponding bilinear alternating map is obtained as an ultraproduct of finite structures. 
It is a well-known example of a supersimple theory of SU-rank 1.
\medskip

We show that the theories of all Fra\"{\i}ss\'e limits $D(n)$ are supersimple of SU-rank 1. To prove this we check the properties of non-forking that characterize simple theories \cite{KP}. Before we show that each group $G$ in ${\mathbb G}(n)$ where $G/Z(G)$ is infinite has the Independence Property especially all $D(n)$.

\section{TP{\boldmath$_2$} of {\boldmath${\rm Th}(D)$}}\label{s1}

\begin{prop}\label{p2.1}
In ${\rm Th}(D)$ the formulae $[x,y_1]=[y_2,y_3]$ has the tree property of the second kind.
\end{prop}

{\em Proof\/}. Since $D$ is the Fra\"{\i}ss\'e limit of ${\mathbb K}_{2,p}^P$ there is an embedding of an infinite free group of $G_{2,p}$ in $D$. Assume that $\{b_\alpha:\alpha<\omega\}\cup\{c_{\alpha,i},d_{\alpha,i}:\alpha<\omega, i<\omega\}$ are free generators of  such an infinite free subgroup. We consider the array
\[
 \ol{a}_{\alpha,i}=\{(b_{\alpha,i}, c_{\alpha,i},d_{\alpha,i}):\alpha<\omega, i<\omega\}
\]
where $b_{\alpha,i}=b_\alpha$ for all $\alpha$ and $i$:
\[
\mbox{Then }\; D\vDash\neg\exists x([x,b_\alpha]=[c_{\alpha,i},d_{\alpha,i}]\wedge[x,b_\alpha]  = [c_{\alpha,j},d_{\alpha,j}])
\]
for fixed $\alpha$ and $i\ne j$. Now let $f$ be any map of $\omega$ into $\omega$. Then the set 
\[
 \{[x,b_\alpha]=[c_{\alpha,f(\alpha)},d_{\alpha,f(\alpha)}]:\alpha<\omega\}
\]
is consistent, since $D$ is a Fra\"{\i}ss\'e limit.\hfill$\Box$

\section{The amalgamation property for {\boldmath${\mathbb K}(n)$}}\label{s3}

Let $G$ be a group in ${\mathbb G}(n)$ with the elements $c_1^G,\ldots,c_n^G$ short $c_1,\ldots,c_n$. Let $P(G)$ be the subgroup generated by $c_1,\ldots,c_n$. In the language $L(n)$ \, $P(G)$ is the $L(n)$-substructure generated by the empty set. By definition $G'\subseteq P(G)\subseteq Z(G)$ and the linear dimension ${\rm ldim}(P(G))$ of $P(G)$ is $n$.
\medskip

In \cite{Bau} a functor $F$ from ${\mathbb G}_{2,p}^P$ into the category ${\mathbb B}^P$ of bilinear alternating maps $(V,W,\beta)$ is defined where $V$ and $W$ are ${\mathbb F}_p$-vector spaces and $\beta$ is a bilinear alternating map from $V$ into $W$. Morphisms of ${\mathbb B}^P$ from $(V_1,W_1,\beta_1)$ to $(V_2,W_2,\beta_2)$ consists of vector space embeddings $f:V_1$ into $V_2$ and $g:W_1$ into $W_2$ that commute with the bilinear maps $\beta_i$.
\[
\begin{xy}
\xymatrix{
V_1\ar[d]^{\textstyle f}&\times&V_1\ar[d]^{\textstyle f}\ar[r]^{\textstyle\beta_1} &W_1\ar[d]^{\textstyle g, }\\
V_2&\times&V_2 \ar[r]^{\textstyle\beta_2}&W_2
}
\end{xy}
\]
$F$ is defined in the following way: $F(G)$ is $(V,W,\beta)$ where $V=G/P(G)$, $W= P(G)$ and $\beta$ is induced by $[\;,\;]$. If $f: G\to H$ then $F(f)=(\ol{f}, f\restriction P)$ where $\ol{f}:G/P(G)\to H/P(H)$ is induced by $f$. 
$F$ is a bijection on the level of objects up to isomorphisms.
\medskip

If we consider the category ${\mathbb G}(n)$, then the morphisms $f:G\to H$ send $c_i^G$ to $c_i^H$. Hence $f$ induces an isomorphism of $P(G)$ onto $P(H)$. We call ${\mathbb B}(n)$ the corresponding category of bilinear alternating maps $(V,P,\beta)$ where $P=\langle c_1,\ldots,c_n\rangle$ is fixed. The morphisms have the form $(g,{\rm id})$. We define the functor $F$ from ${\mathbb G}(n)$ to ${\mathbb B}(n)$ as above and obtain as in \cite{Bau2}:

\begin{lemma}\label{l3.1}
\begin{enumerate}
 \item[{\rm i)}] $F$ is a functor of ${\mathbb G}(n)$ onto ${\mathbb B}(n)$ that is a bijection for the objects of the categories up to isomorphisms.

\item[{\rm ii}] If $G_0 \in {\mathbb G}(n)$ and $(g,id)$ is ann embedding of $F(G_0)$ into some $(V,P,\beta)$, then there are some $G \in {\mathbb G}(n)$ and some embedding $f$ of $G_0$ into $G$, such that $F(G) = (V,P, \beta)$ and $F(f) = (g,id)$.

\item[{\rm iii)}] In ${\mathbb G}(n)$ we consider $e_0:G_0\to G$, $e_1:H_0\to H$ where $f_0$ is an isomorphism of $G_0$ onto $H_0$. In ${\mathbb B}(n)$ we assume that there is $g$ such that
\[
\begin{xy}
\xymatrix{
F(G_0)\ar[d]^{F(f_0)}\ar[r]^{F(e_0)} &F(G)\ar[d]^{(g,{\rm id})}\\
F(H_0) \ar[r]^{F(e_1)}&F(H)\;\;.
}
\end{xy}
\]
\end{enumerate}
Then there is an embedding $f$ of $G$ into $H$ such that $F(f)=(g,{\rm id})$ and 
\[
\begin{xy}
\xymatrix{
G_0\ar[d]^{f_0}\ar[r]^{e_0} &G\ar[d]^{f}\\
H_0 \ar[r]^{e_1}&H\;\;.
}
\end{xy}
\]
\end{lemma}

Lemma~\ref{l3.1} shows that AP for ${\mathbb B}(n)$ implies AP for ${\mathbb G}(n)$ as in \cite{Bau}. To show AP for ${\mathbb B}(n)$ we cannot use the free amalgam as in \cite{Bau}.
\smallskip

Assume
\begin{eqnarray*}
(f_A,{\rm id}):&&(V_B,P,\beta_B)\longrightarrow(V_A,P,\beta_A),\\
(f_C,{\rm id}):&&(V_B,P,\beta_B)\longrightarrow(V_C,P,\beta_C).
\end{eqnarray*}
W.l.o.g. $V_B$ is a common subspace of $V_A$ and $V_C$. Let $V_D$ be the vector space amalgam 
$V_C\bigoplus\limits_{V_B} V_A$  with respect to $f_A$ and $f_C$. We get the desired amalgam $\langle V_D,P,\beta_D\rangle$ if 
\[
\beta_D=\beta_A\;\mbox{ on }\; V_A \quad\mbox{ and } \quad \beta_D=\beta_C\;\mbox{ on }\;V_C
\]
and the rest is obtained in the following way: 
If $X$ is a basis of $V_A$ over $V_B$ and $Y$ is a basis of $V_C$ over $V_B$ then we can choose for each pair $x\in X$ and $y\in Y$ \, $\beta_0(x,y)$ in $P$ as we want.
\medskip

In our context AP implies JEP.

\begin{theorem}\label{t3.1}
${\mathbb K}(n)$ has HP, JEP and AP. Hence the Fra\"{\i}ss\'e limit $D(n)$ exists. It is $\aleph_0$-categorical. ${\rm Th}(D(n))$ has the elimination of quantifiers.
\end{theorem}
The theorem uses the known theory. See \cite{Ho}. Uniform local finiteness and finite signature for ${\mathbb K}(n)$ imply $\aleph_0$-categoricity and elimination of quantifiers. ${\rm Th}(D(n))$ can be axiomatized by the following sentences: Let $M$ be a model of ${\rm Th}(D(n))$.
\begin{enumerate}
 \item[$\Sigma\,1)$] $M$ is a nilpotent group of class 2 with exponent $p$.
\item[$\Sigma\,2)$] $M'=Z(M)=\langle c_1,\ldots,c_n\rangle$ is of linear dimension $n$. 
\item[$\Sigma\,3)$] For $B\subseteq A$ in ${\mathbb K}(n)$ it holds: If $B'\subseteq M$ and $B'\cong B$, then this embedding of $B$ into $M$ can be extended to $A$.
\end{enumerate}

In the case $n=1$ these axioms imply that $M$ is infinite and $M'=Z(M)$ is cyclic.
By U.~Felgner \cite{Fe} $D(1)$ is the extra special $p$-group, since his axiomatization is $\Sigma\,1)$ $M'=Z(M)$ is cyclic and infiniteness.
\bigskip

{\bf Question} \, Is there an easier axiomatization of ${\rm Th}(D(n))$ for $n\ge 2$?

\section{Independence property in {\boldmath${\mathbb G}(n)$}}\label{s4}

Assume $M\vDash {\rm Th}(D(1))$ and $M$ is countable. We write $c$ instead of $c_1$. By \cite{Fe} $M$ is a central product over $\langle c\rangle$:
\[
 M=\bigodot\limits_{\stackrel{\scriptstyle\langle c\rangle}{i<\omega}}\langle c,a_ib_i\rangle
\]
where $c$ is a generator of the cyclic subgroup $M'=Z(M)$ and $[b_i,a_i]=c$.
\smallskip

By the elimination of quantifiers of ${\rm Th}(D(1))$ \, $a_0\hat{\;}b_0,a_1\hat{\;}b_1\ldots,a_n\hat{\;}b_n,\ldots$ is an indiscernible sequence in $M$. Then 
$a_1\hat{\;}b_1,a_2\hat{\;}(b_2\circ b_0),a_3\hat{\;}b_3,a_4\hat{\;}(b_4\circ b_0),\ldots$ and $b_1,b_2\circ b_0,b_3,b_4\circ b_0,\ldots$ are indiscernible sequences in $M$.
We have $M\vDash[b_{2i+1},a_0]=1$ for $i<\omega$ and $M\vDash[b_{2i}\circ b_0,a_0]=c$ for $i\le i<\omega$. 
\medskip

We have shown  (see \cite{A}):

\begin{lemma}\label{l4.1}
The formula $[y,x]=1$ has the independence property in ${\rm Th}(D(1))$.
\end{lemma}

Let $G$ be in ${\mathbb G}(n)$ with $G/Z(G)$ is infinite. For $a\in G\setminus Z(G)$ choose a maximal linearly independent subset $\{e_1,\ldots,e_m\}=X_a$ of $P(G)$ such that for every $1\le i\le m\le n$ there is some $b_i\in G$ with $[a,b_i]=e_i$. Let $E_a$ be $\{b_1,\ldots,b_m\}$. If $[a,b]=t\ne 1\in P(G)$, then $t=\sum\limits_{1\le i\le m} e_i^{r_i}$ and $[a,b\cdot b_1^{p-r_1}\cdot\ldots\cdot b_m^{p-r_m}]=1$. Hence every element $a\in G$ has a centralizer $C(a)$ of index $\le n$ and $G=\langle a,E_a\rangle\circ C(a)$. 
\medskip

Now we start again with $d_0\in G\setminus Z(G)$. Then $X_{d_0}\ne \emptyset$ and $E_0=E_{d_0}\ne\emptyset$ and we choose $e_0\in E_0$ with $[d_0,e_0]\ne 1$.
Since $C(\langle d_0,E_0\rangle)$ has finite index in $G$ there is some
\[
 d_1\in C(\langle d_0,E_0\rangle)\quad\mbox{ with }\; d_1\not\in Z(G).
\]
We get $E_1=E_{d_1}\ne\emptyset$ and choose $b_1\in E_1$. 
We can repeat this argument and get 
\[
 d_2\in C(\langle d_0,e_0,E_0,E_1\rangle), \quad d_2\not\in Z(G).
\]
Finally we have $d_0,e_0, d_1,e_1,\ldots$ with $[d_i,e_i]\ne 1$ and $[d_i,d_j]=1$, $[d_i,e_j]=1$ and $[e_i,e_j]=1$ for $i\ne j$.
We can select a subsequence with $[d_i,e_i]=c\ne 1$ for some $c\in P(G)$. Assume w.l.o.g. $[d_i,e_i]=c$ for all $i<\omega$. We have shown that $D(1)$ is a subgroup of $G$. Since the independence property of $D(1)$ is given by a quantifier formula we get

\begin{theorem}\label{t4.2}
For every $G\in{\mathbb G}(n)$ with $G/Z(G)$ infinite we have:
\begin{enumerate}
 \item[{\rm i)}] There is an embedding of $D(1)$ in $G$.
\item[{\rm ii)}] $G$ has the independence property.
\end{enumerate}
\end{theorem}

\begin{corollary}\label{c4.3}
The Fra\"{\i}ss\'e limits $D(n)$ of ${\mathbb K}(n)$ have the independence property.
\end{corollary}

\section{Superstability of {\boldmath$Th(D(n))$}}\label{s5}

Let ${\mathbb C}(n)$ be a monster model of $Th(D(n))$. We define
\[
A\ind\limits_B{}\!\!^0 C, \mbox{ if }\;  \langle A\rangle\cap\langle C\rangle=\langle B\rangle.
\]
Note that all substructures as $\langle A\rangle$ contain $P({\mathbb C}(n))$.
\smallskip

We have to check that $\ind^0$ fulfils the conditions of B.~Kim and A.~Pillay \cite{KP} that characterize Non-forking. Working in the vector space ${\mathbb C}(n)/P({\mathbb C}(n))$ Monotonicity, Transitivity, Symmetry, Finite Character, and Local Character are easily shown.
\medskip

{\bf Existence:} $\ol{a},B\subseteq A$ are considered in ${\mathbb C}(n)$. Then there is some $\ol{d}$ in ${\mathbb C}(n)$ with ${\rm tp}(\ol{a}/B)={\rm tp}(\ol{d}/B)$ and $\ol{d}\ind\limits_B{}\!\!^0A$.
\medskip

W.l.o.g. $B$ and $A$ are $L(n)$-substructures. Since $P\subseteq B$ we can assume that $\ol{a}$ is linearly independent over $B$. Choose $X_B$ and $X_A$ such that the images of $X_B$ and $X_BX_A$ are vector space bases of $B/P$ and $A/P$, respectively.
Let $\beta((X_BX_A)^2)$ be the set of all $\beta(b_1,b_2)=[b_1,b_2]$ where $b_1,b_2\in X_BX_A$. Then $A$ is uniquely determined by $X_BX_A$ and $\beta((X_BX_A)^2)$.
\smallskip

Now we define an extension $G$ of $A$. Let $\ol{e}$ linearly independent over $A$. $\ol{e}X_BX_A$ is linearly independent over $P$. $\beta((\ol{e}\,\hat{\;}X_BX_A)^2)$ is chosen as any extension of $\beta((X_BX_A)^2)$ and $\beta((\ol{e}X_B)^2)$, where the last set is obtained from $\beta((\ol{a}X_B)^2)$ by replacing $a_i$ in $\ol{a}$ by $e_i$ in $\ol{e}$.  $G=\langle\ol{e}A\rangle$ is a structure in ${\mathbb K}(n)$. By the axioms $\Sigma\,3)$ of ${\rm Th}(D(n))$ there is an embedding of $\ol{e}$ onto $\ol{d}$ over $A$ in ${\mathbb C}(n)$ . By quantifier elimination ${\rm tp}(\ol{d}/B)={\rm tp}(\ol{a}/B)$. Furthermore $\ol{d}\ind\limits_B{}\!\!^0A$ by construction. 
\medskip

Finally we have to show:

\subsubsection*{Independence over Models}

Let $M\preceq {\mathbb C}(n)$, ${\rm tp}(\ol{a}{}^0/M)={\rm tp}(\ol{a} {}^1/M)$
\[
\ol{b}{}^0\ind\limits_M{}\!\!^0\;\ol{b}{}^1,\qquad \ol{a}{}^0\ind\limits_M{}\!\!^0\;\ol{b}{}^0,\qquad
\ol{a}{}^1\ind\limits_M{}\!\!^0\;\ol{b}{}^1.
\]
Then there is some $\ol{e}$ with 
\[
 {\rm tp}(\ol{e}/M\ol{b}{}^0)={\rm tp}(\ol{a}_0/M\ol{b}{}^0),\qquad {\rm tp}(\ol{e}/M\ol{b}{}^1)={\rm tp}(\ol{a}_1/M\ol{b}{}^1),
\]
and
\[
 \ol{e}\ind\limits_M{}\!\!^0\;\ol{b}{}^0\ol{b}{}^1.
\]
Let $X_M$ be a set in $M$ such that its image is a vector space basis of $M/P$. By assumption we can assume that w.l.o.g.  $\ol{b}{}^0\ol{b}{}^1$ is linearly independent over $M$ modulo~$P$. We choose $\ol{d}$ linearly independent over $\ol{b}{}^0\ol{b}{}^1X_M$ modulo~$P$. Now we extend $\langle \ol{b}{}^0,\ol{b}{}^1,X_M\rangle$ to a group $G$ in ${\mathbb G}(n)$ defined on $\ol{d}\,\ol{b}{}^0\ol{b}{}^1X_M$. We extend $\beta((\ol{b}{}^0\ol{b}{}^1X_M)^2)$ to $\beta((\ol{d}\,\ol{b}{}^0\ol{b}{}^1X_M)^2)$ by the following: 
\begin{eqnarray*}
\beta(d_i,m)&&\mbox{for }\:d_i\in\ol{d}\;\mbox{ and }\; m\in X_M\;\mbox{ is given by }\;\beta(a_i^0,m)=\beta(a_i^1,m),\\
\beta(d_i,b_j^0)&&\mbox{is given by }\;\beta(a_i^0,b_j^0), \mbox{ and}\\
\beta(d_i,b_j^1)&&\mbox{is given by }\; \beta(a_i^1,b_j^1).
\end{eqnarray*}
Now we find an image $\ol{e}$ of $\ol{d}$ in ${\mathbb C}(n)$ over $\langle\ol{b}{}^0,\ol{b}{}^1,M\rangle$ by axioms $\Sigma\,3)$ that defines an embedding. By elimination of quantifiers and the consltruction $\ol{e}$ has the desired properties.

\begin {theorem}\label{t5.1}
$\ind^0$ is non-forking for $D(n)$. $D(n)$ is supersimple of SU-rank $1$. It is not stable.
\end {theorem}

{\em Proof\/}. As shown above $\ind^0$ is non-forking and $D(n)$ is simple. Any type ${\rm tp}(\ol{a}/A)$ does not fork on a finite subset of $A$. By the description of non-forking we have SU-rank 1. In Chapter~4 it is shown that $D(n)$ has the independence property.

\end{document}